\documentclass{amsart}

\usepackage{amssymb}

\usepackage{tikz-cd}

\usepackage{mathrsfs}

\usepackage{float}

\usepackage{mathtools}
\DeclarePairedDelimiter{\floor}{\lfloor}{\rfloor}

\usepackage[noabbrev,capitalise]{cleveref}

\usepackage{amsrefs}

 \newcommand{\sE}{\mathscr{E}}
\newcommand{\sF}{\mathscr{F}}
\newcommand{\sG}{\mathscr{G}}
\newcommand{\sI}{\mathscr{I}}
\newcommand{\sO}{\mathscr{O}}
\newcommand{\sP}{\mathscr{P}}

\newcommand{\calF}{\mathcal{F}}

\newcommand{\All}{\mathscr{A}\ell\ell}

\newcommand{\etilde}{\tilde{e}}

\DeclareMathOperator{\res}{res}
\DeclareMathOperator{\Id}{Id}
\DeclareMathOperator{\Res}{Res}
\DeclareMathOperator{\Ind}{Ind}
\DeclareMathOperator{\Sp}{Sp}
\DeclareMathOperator{\sk}{sk}
\renewcommand{\Im}{\operatorname{Im}}

\newcommand{\Proj}{\mathrm{P}}

\newcommand{\op}{\mathrm{op}}
\newcommand{\sw}{\mathbb{D}}

\newcommand{\pt}{\mathrm{pt}}

\newcommand{\Integers}{\mathbf{Z}}
\newcommand{\Reals}{\mathbf{R}}
\newcommand{\Complex}{\mathbf{C}}
\newcommand{\HH}{\mathbf{H}}
\newcommand{\Ffield}{\mathbf{F}}

\DeclareMathOperator*{\hocolim}{hocolim}
 \theoremstyle{plain}
\newtheorem{Theorem}{Theorem}[section]
\newtheorem{Bonus Theorem}[Theorem]{Bonus Theorem}

\newtheorem{Lemma}[Theorem]{Lemma}

\newtheorem{Proposition}[Theorem]{Proposition}

\newtheorem{Corollary}[Theorem]{Corollary}

\newtheorem{Unverified Claim}[Theorem]{Unverified Claim}
\newtheorem{TheoremA}{Theorem}

\Crefname{TheoremA}{Theorem}{Theorems}
\theoremstyle{definition}
\newtheorem{Definition}[Theorem]{Definition}

\theoremstyle{remark}

\newtheorem{Example}[Theorem]{Example}

\newtheorem{Non-example}[Theorem]{Non-example}
\newtheorem{Remark}[Theorem]{Remark}

\newtheorem{Notation}[Theorem]{Notation}

\begin{document}
\title{Quantifying Quillen's uniform $\calF_p$-isomorphism theorem}
\author{Koenraad van Woerden}
\email{koenvanwoerden@gmail.com}
\thanks{The author was partly supported by the SFB 1085 -- Higher Invariants, Regensburg.}
    \begin{abstract}
        Let $G$ be a finite group with $2$-Sylow subgroup of order less than or equal to 16. For such a $G$, we prove a quantified version of Quillen's uniform $\calF_p$-isomorphism theorem, which holds uniformly for all $G$-spaces.

        We do this by bounding from above the exponent of Borel equivariant $\Ffield_2$-cohomology, as introduced by Mathew--Naumann--Noel, with respect to the family of elementary abelian 2-subgroups.
    \end{abstract}

    \maketitle

    \section{Introduction}
    For a finite group $G$, consider a cohomology class $u \in H^*(BG; \Ffield_p)$ that restricts to zero on all elementary abelian $p$-subgroups (i.e.,\ groups of the form $(\Integers/p)^{\times l}$ for some integer $l \geq 0$). It is a theorem of Quillen that $u$ is nilpotent. In fact, Quillen showed
\begin{Theorem}[{{\cite[Thm.\ 6.2]{quillen71spectrum}}}]
    \label{thm:quillenfp}
    \label{thm:quillenborfp}
    For $X$ any paracompact $G$-space of finite cohomological dimension, the map
    \begin{equation}
      \widetilde{\res} \colon H_G^*(X;\Ffield_p) \to \lim_{E\subset G \text{ el.\ ab.\ $p$-gp.}} H_E^*(X;\Ffield_p),
      \label{eq:quillenfp}
\end{equation}
where the maps in the indexing category for the limit are given by restricting along subgroups and conjugation, is a uniform $\calF_p$-isomorphism, which means that there is an integer $n \geq 0$ such that
\begin{enumerate}
    \item Every $u \in \ker(\widetilde{\res})$ satisfies $u^n = 0$.
    \item Every $v \in \lim_E H^*_E(X; \Ffield_p) \setminus \Im(\widetilde{\res})$ satisfies $v^{p^n} \in \Im(\widetilde{\res})$. 
\end{enumerate}

\end{Theorem}
Here $H_G^*(X; \Ffield_p)$ denotes the mod-$p$ Borel equivariant cohomology of a $G$-space $X$.
These results led to many structural results in group cohomology. Quillen himself immediately deduced:
\begin{Corollary}[{{\cite[Cor.\ 7.8]{quillen71spectrum}}}]
    The Krull dimension of $H^*(BG;\Ffield_p)$ equals the rank of the maximal elementary abelian $p$-subgroup of $G$.
\end{Corollary}
Results directly building on \cref{thm:quillenborfp} include a theorem of Duflot on the depth of $H_G^*(X;\Ffield_p)$ \cite[Thm.\ 1]{duflot81}, a theorem on the complexity of $kG$-modules by Alperin--Evens \cite{alperinevens81}, Benson's description of the image of the transfer map \cite[Thm.\ 1.1]{benson93}, and a theorem on the depth of group cohomology rings by Carlson \cite[Thm.\ 2.3]{carlson95}. These results, and Quillen's original result all indicate the importance of the elementary abelian $p$-subgroups in Borel equivariant $\Ffield_p$-cohomology in general, and group cohomology with $\Ffield_p$-coefficients in particular.

It is natural to ask what one can say about the $n$ in \cref{thm:quillenfp}. One approach to this question is to apply the work of Kuhn \cites{kuhn07, kuhn13}, which builds on work of Henn--Lannes--Schwartz \cite{hls95}. In particular, for $p=2$, this gives explicit upper bounds on the nilpotence degree of the kernel of (\ref{eq:quillenfp}), for $X = \pt$. Alternatively, for the case $X = \pt$, one can explicitly compute both sides of (\ref{eq:quillenfp}) to determine a feasible $n$.

The case $X \neq \pt$ is not covered by these bounds however. We will consider all $X$, not necessarily $X = \pt$, by using the approach introduced by Mathew--Naumann--Noel in \cite{mnn}, where the map (\ref{eq:quillenfp}) is realized as the edge homomorphism of a homotopy limit spectral sequence, called the $\sF$-homotopy limit spectral sequence
\begin{equation*}
    E_2^{s,t} = \sideset{}{^s} \lim_{\sO(G)_{\sF}^{\op}} H_H^t(X; \Ffield_p) \Rightarrow H_G^*(X;\Ffield_p)
\end{equation*}
converging strongly to the target \cite[Prop.\ 2.24]{mnn}, where $X$ can be any $G$-space (or more generally, a $G$-spectrum), and $\sF$ any family of subgroups of $G$ which contains at least the family $\sE_{(p)}$ of elementary abelian $p$-subgroups. The indexing category is the subcategory of the orbit category $\sO(G)$ spanned by the orbits $G/H$ with $H$ in the family $\sF$. The key property of the $\sF$-homotopy limit spectral sequence is that it collapses at a finite page with a horizontal vanishing line \cite[Thm.\ 2.25]{mnn}. 

This implies that every computation with the $\sF$-homotopy limit spectral sequence is a finite one. Moreover, in many concrete situations we can establish a bound on the height of the horizontal vanishing line, and a bound on which page it appears. Besides implying when the $\sF$-homotopy limit spectral sequence will have collapsed, it can also be used to deduce differentials. This was illustrated in the computation of the cohomology of the quaternion group of order 8 in \cite[Ex.\ 5.18]{mnn}. The  computational utility of the $\sF$-homotopy limit spectral sequence will be illustrated in a forthcoming paper, by using it to compute the cohomology of all 2-groups up to order 16. 

Varying $X$ over all $G$-spectra, this horizontal vanishing line turns out to have a uniform upper bound in height \cite[Prop.\ 2.26]{mnn}. The minimal upper bound of this height is one of the definitions of the $\sE_{(p)}$-exponent $\exp_{\sE_{(p)}}\underline{H\Ffield_p}_G$. An equivalent definition is the following.
\begin{Definition}[{{\cite[Prop.\ 2.26]{mnn}}}]
    The $\sE_{(p)}$-exponent of $\underline{H\Ffield_p}$ is the minimal integer $n \geq 0$ such that there exists an $n$-dimensional CW-complex $X$ with isotropy in $\sE_{(p)}$ such that the canonical map $\underline{H\Ffield_p} \to F(X_+, \underline{H\Ffield_p})$ admits a retraction
    \begin{equation*}
        \underline{H\Ffield_p} \to F(X_+, \underline{H\Ffield_p}) \to \underline{H\Ffield_p}.
    \end{equation*}
\end{Definition}
In practice one can often determine this $\sE_{(p)}$-exponent, and this leads to a quantified version of \cref{thm:quillenfp}, because one has $n \leq \exp_{\sE_{(p)}}\underline{H\Ffield_p}_G$ \cite[Thm.\ 3.20, Rem.\ 3.22]{mnn}. The identification of $\sE_{(2)}$-exponents for 2-groups is the principal goal of this paper, and leads to the main theorem:
\begin{TheoremA}
    \label{thm:A}
    Let $G$ be a finite group with a 2-Sylow subgroup of order $\leq 16$, $X$ any $G$-space, and let $I$ be the kernel of
    \begin{equation*}
            \widetilde{\res} \colon H_G^*(X;\Ffield_2) \to \lim_{E\subset G \text{ el.\ ab.\ $2$-gp.}} H_E^*(X;\Ffield_2).
    \end{equation*}
    Then $I^4 = 0$. Moreover, if $u$ is any element in the codomain of $\widetilde{\res}$, then $u^{8}$ is in the image of $\widetilde{\res}$.
\end{TheoremA}
\Cref{thm:A} follows from combining \cite[Thm.\ 3.20]{mnn} and \cref{lem:psylowhfpexp} with the upper bounds on the exponents from the following theorem:
\begin{TheoremA}
    The exponents of Borel equivariant $\Ffield_2$-cohomology for the groups of order less than or equal to 16 are bounded above by the values in the following table:
\begin{table}[H]
    \centering
     \begin{tabular}{r|r|r}
         $G$ & $\exp_{\sE_{(2)}} \underline{H\Ffield_2}_G$ & Reference \\
         $e$ & $1$ & \cref{prop:expelab} \\
         \hline
         $C_2$ & $1$ & \cref{prop:expelab} \\
         \hline
         $C_2 \times C_2$ & $1$ & \cref{prop:expelab} \\
         $C_4$ & $2$ & \cref{prop:expelab} \\
        $C_2^{\times 3}$ & $1$ & \cref{prop:expelab} \\
        $C_2 \times C_4$ & $2$ & \cref{prop:expelab} \\
        $C_8$ & $2$ & \cref{prop:expelab} \\
        $D_8$ & $2$ & \cref{cor:expd8upperbound} \\
        $Q_8$ & $4$ & \cite[Ex.\ 5.18]{mnn} \\
        \hline
        $C_2^{\times 4}$ & $1$ & \cref{prop:expelab}  \\
        $C_2^{\times 2} \times C_4$ & $2$ & \cref{prop:expelab}  \\
        $C_4 \times C_4$ & $3$ & \cref{prop:expelab} \\
        $C_8 \times C_2$ & $2$ & \cref{prop:expelab}  \\
        $C_{16}$ & $2$ & \cref{prop:expelab}   \\
        $D_{16}$ & $2$ & \cref{cor:expd8upperbound}\\
        $Q_{16}$ & $4$ & \cref{q2n:expupper}  \\
        $ SD_{16} =  C_8 \overset{3}{\rtimes} C_2$ & $4$ & \cref{prop:sd16upper}  \\
        $ M_{16} = C_8 \overset{5}{\rtimes} C_2$ & $4$ & \cref{m16:expupper}  \\
        $D_8 \ast C_4$ & $4$ & \cref{d8cpc4:elexpupperbound2} \\
        $C_4 \rtimes C_4$ & $4$ & \cref{prop:c4xc2sdc2exp} \\
        $(C_4 \times C_2) \overset{\psi_5}{\rtimes} C_2$ & $2$ & \cref{prop:c4xc2sdc2exp} \\
        $Q_8 \times C_2$ & $4$ & \cref{prop:q8xc2e2exp} \\
        $D_8 \times C_2$ & $2$ & \cref{prop:d8xc2e2exp} \\
       \end{tabular}
    \caption{Upper bounds on the $\sE_{(2)}$-exponents of the groups of order $\leq 16$.}
    \label{table:groupsleq16}
\end{table}
\label{thm:B}
\end{TheoremA}
For a description of the groups appearing in \cref{thm:B} we refer to the sections of the referred propositions. The first column lists the groups of order $\leq 16$, the second column an upper bound on the exponent, and the last column gives a forward reference for the claim.
\begin{Remark}
    All these upper bounds are in fact equalities, except possibly when $G$ equals $SD_{16}$, $M_{16}$ or $C_4 \rtimes C_4$, in which the exponent could be 3. These lower bounds will be part of the content of the forthcoming paper mentioned above.
\end{Remark}
For specific 2-Sylow subgroups of order less than or equal to $16$, one can obtain an improved version of \cref{thm:A} by using the upper bound on the relevant exponent from \cref{thm:B}.

\subsection{Organization}
The proofs of the main results are contained in \cref{sec:expsm2gps}. \cref{sec:fnilp} summarizes what we need from \cites{mnn, mnnnd}, and \cref{ch:explems} contains lemmas that are used to prove the main results.

\subsection{Acknowledgements}
I would like to thank Akhil Mathew for a detailed explanation of the use of the projective bundle theorem in \cite[Ex.\ 5.18]{mnn}.

This paper is based on the author's PhD thesis, which was written under the supervision of Justin Noel. I would like to thank him for his guidance and support, and for his help with the statement and proofs of the exponent lemmas in \cref{ch:explems}.

Finally, I would like to thank the referee for providing helpful remarks on this paper.
\subsection{Notation and conventions}
Throughout $G$ denotes a finite group. If a non-equivariant cohomology theory is represented by a spectrum $E$, then we denote the spectrum representing the Borel equivariant version of this cohomology theory by $\underline{E}_G$. In particular, we denote Borel-equivariant $\Ffield_2$-cohomology $\underline{H\Ffield_2}_G$. We may omit the $G$ from the notation if the group $G$ is clear from the context. For a finite group $G$ and a prime $p$, we denote the family of elementary abelian $p$-subgroups of $G$ by $\sE(G)_{(p)}$ or $\sE_{(p)}$ if the group is clear from the context.

     \section{$\sF$-nilpotence}
    \label{sec:fnilp}
    We recall the notion of $\sF$-nilpotence from \cite{mnn}, for which we first need to recall the following space.
\begin{Definition}[{{\cite[Satz 1 and Beweis, Bemerkung 2]{tomDieck72}}}]
    Let $\sF$ be a family of subgroups of $G$. Then the \textbf{universal $\sF$-space} $E\sF$ is the $G$-space given by the homotopy colimit
    \begin{equation*}
        E\sF = \hocolim_{\sO(G)_{\sF}} G/H.
    \end{equation*}
\end{Definition}
The space $E\sF$ is, up to $G$-equivalence, characterized by
\begin{equation*}
    E\sF^H \simeq \begin{cases}
        \varnothing & \text{if $H \not\in \sF$,} \\
        \pt & \text{if $H \in \sF$.}
    \end{cases}
\end{equation*}
Using the space $E\sF$, we can give one of the equivalent definitions of $\sF$-nilpotence.
\begin{Definition}[{{Cf. \cite[Def.\ 1.4]{mnn}}}]
    Let $M$ be a $G$-spectrum. Then $M$ is said to be $\sF$-nilpotent if there is an $n$ such that $M$ is a retract of $F(\sk_{n-1}E\sF_+,M)$. The minimal $n \geq 0$ for which this holds is called the $\sF$-exponent of $M$, and denoted $\exp_{\sF}M$. 
\end{Definition}
Being $\sF$-nilpotent is a strong condition, it implies for example the following.
\begin{Proposition}[{{\cite[Prop.\ 2.8]{mnn}}}]
    Let $M$ be $\sF$-nilpotent. Then $M$ is $\sF$-complete and $\sF$-colocal, that is, the $\sF$-completion map $M \to M(E\sF_+,M)$ and the colocalization map $M \wedge E\sF_+ \to M$ are weak equivalences.
\end{Proposition}
A trivial example is the following.
\begin{Example}
    Let $\All$ be the family of all subgroups of a group $G$. Every $G$-spectrum $M$ is $\All$-nilpotent with $\exp_{\All} M \leq 1$, and $\exp_{\All} M = 0$ if and only if $M$ is contractible.
    \label{exam:allnilp}
    \label{prop:allexp}
    \label{lem:expallsubgroups}
\end{Example}
The following two propositions are immediate from the above definition.
\begin{Proposition}[{{\cite[Prop. 6.39]{mnnnd}}}]
    A $G$-spectrum $M$ is $\sF$-nilpotent and $\sG$-nilpotent if and only if $M$ is $\sF \cap \sG$-nilpotent.
    \label{prop:intersectnilp}
\end{Proposition}
\begin{Proposition}[{{\cite[Def.\ 1.4]{mnn}}}]
    If $M$ is $\sF$-nilpotent and $\sG \supset \sF$, then $M$ is $\sG$-nilpotent.
    \label{prop:supfamnilp}
\end{Proposition}
Combining \cref{exam:allnilp}, \cref{prop:intersectnilp} and \cref{prop:supfamnilp} shows that every $G$-spectrum $M$ has a minimal family $\sF$ such that $M$ is $\sF$-nilpotent (cf.\ the remark after \cite[Def. 1.4]{mnn}). 
\begin{Definition}[{{\cite[remark after Def.\ 1.4]{mnn}}}]
    The minimal family $\sF$ such that a $G$-spectrum $M$ is $\sF$-nilpotent is called the \textbf{derived defect base} of $M$.
\end{Definition}
The following is the main case of interest for us.
\begin{Proposition}[{{\cite[Prop.\ 5.16]{mnn}}}]
    For $G$ any finite group, the derived defect base of $\underline{H\Ffield_2}_G$ is $\sE_{(2)}$.
    \label{prop:derdefbaseborel}
\end{Proposition}
We can now state the uniform upper bound on the height of the horizontal vanishing line of the $\sF$-homotopy limit spectral sequence and the page on which it appears.
\begin{Proposition}[{{\cite[Prop.\ 2.26, Rem.\ 2.27]{mnn}}}]
    \label{expprop}
    Let $G$ be a finite group, $\sF$ a family of subgroups, and $M$ an $\sF$-nilpotent $G$-spectrum. Then the following integers are equal:
    \begin{enumerate}
        \item The $\sF$-exponent of $M$.
        \item The minimal $N$ such that for all $G$-spectra $X$, the $\sF$-homotopy limit spectral sequence $E_*^{*,*}(X)$ admits a vanishing line of height $N$ on the $E_{N+1}$-page: $E_{N+1}^{s,*} = E_\infty^{s,*} = 0$ for all $s \geq N$.
        \item The minimal $n$ such that the canonical map $F(E\sF_+,M) \simeq M \to F(\sk_{n-1}E\sF_+, M)$ admits a retraction.
        \item The minimal $n'$ such that there is an $(n'-1)$-dimensional CW-complex $X$ with isotropy in $\sF$ such that $M$ is a retract of $F(X_+, M)$.
        \item The minimal $m$ such that the canonical map $\sk_{m-1}E\sF \wedge M \to M$ admits a section. \label{expprop:smash2}
        \item The minimal $m'$ such that there is an $(m'-1)$-dimensional CW-complex $X$ with isotropy in $\sF$ such that $M$ is a retract of $X_+ \wedge M$. \label{expprop:smash}
    \end{enumerate}
    Moreover, if $M'$ is any $G$-spectrum, then the existence of an integer for $M'$ as in any one of the items from (2) to (6) implies that $M'$ is $\sF$-nilpotent.
\end{Proposition}
We end this section by recalling from \cite{mnn} some properties of exponents that will be used in the next chapter to prove lemmas about exponents.
\begin{Proposition}
    Let $H \in \sF$. Then $G/H_+$ is $\sF$-nilpotent with $\exp_{\sF} G/H_+ = 1$.  \label{prop:fingsetexp} \label{lem:expgsets}
\end{Proposition}
\begin{Proposition}[{{\cite[Cor.\ 4.15]{mnnnd}}}]
    If $M$ is an $\sF$-nilpotent spectrum and $X$ is any $G$-spectrum, then $F(X,M)$ is $\sF$-nilpotent with $\exp_{\sF} F(X,M) \leq \exp_{\sF} M$.
\end{Proposition}
\begin{Proposition}
    If $N$ is an $\sF$-nilpotent $G$-spectrum and $M$ is any $G$-spectrum, then $M \wedge N$ is $\sF$-nilpotent with $\exp_{\sF}M \wedge N \leq \exp_{\sF} N$.
    \label{prop:smashexp}
\end{Proposition}
\begin{Proposition}[{{\cite[Prop.\ 4.9]{mnnnd}}}]
 \mbox{}
    \begin{enumerate}
        \item If $M$ is a retract of an $\sF$-nilpotent spectrum $N$, then $M$ is $\sF$-nilpotent and $\exp_{\sF}M \leq \exp_{\sF} N$.
        \item If $M'$ and $M''$ are $\sF$-nilpotent and $M' \to M \to M''$ is a cofiber sequence, then $M$ is $\sF$-nilpotent and $\exp_{\sF} M \leq \exp_{\sF} M' + \exp_{\sF} M''$.
    \end{enumerate}
    \label{prop:retrexp}
    \label{lem:expretracts}
\end{Proposition}
\begin{Proposition}
    Let $M_{\alpha}$ be a set of $\sF$-nilpotent spectra with $\sF$-exponents bounded uniformly by $n$. Then $\bigvee_{\alpha} M_\alpha$ is $\sF$-nilpotent with $\sF$-exponent $\leq n$.
    \label{prop:wedgeexp}
\end{Proposition}
\begin{Proposition}
    Let $X$ be an $(n-1)$-dimensional $G$-CW-spectrum with isotropy in $\sF$. Then $X$ is $\sF$-nilpotent and $\exp_{\sF} X \leq n$.
\end{Proposition}
\begin{proof}
    Use induction and the previous propositions.
\end{proof}
\begin{Proposition}
    Let $X$ be a finite dimensional $G$-CW-spectrum with isotropy in $\sF$. Then the equivariant Spanier--Whitehead dual $\sw(X)$ of $X$ is $\sF$-nilpotent, and $\exp_{\sF}\sw(X) = \exp_{\sF} X$.
    \label{prop:cwexp}
    \label{lem:cwexp}
\end{Proposition}
\begin{proof}
    Write $n = \exp_{\sF} X$.
    Let $Y_+$ be an $(n-1)$-dimensional $G$-CW complex with isotropy in $\sF$ such that there is a retraction
    \begin{equation*}
        X \to F(Y_+, X) \simeq \sw(Y_+) \wedge X \to X.
    \end{equation*}
    Applying $\sw(-)$ to this retraction exhibits $\sw(X)$ as a retraction of $Y_+ \wedge \sw(X)$, which has exponent $\leq n$ by \cref{prop:cwexp} and \cref{prop:smashexp}. Therefore, $\exp_{\sF}\sw(X) \leq \exp_{\sF}X$, and replacing $X$ by $\sw(X)$ in this inequality shows equality.
\end{proof}

     \section{Exponent lemmas}
        \label{ch:explems}
        We discuss some lemmas that will be of use in determining exponents. Some of these statements appear as exercises in \cite[Sec.\ 4]{mnnnd}.
         \subsection{Lemmas for $\sF$-exponents}
        \label{sec:exponentlemmas}
We now give various lemmas which describe how $\sF$-exponents can change as the family $\sF$ varies. We have the following two basic product and restriction formulas.
\begin{Remark}
    Let $\sF_1$, $\sF_2$ be two families of subgroups. Then $E(\sF_1 \cap \sF_2)_+ \simeq E\sF_{1_+} \wedge E\sF_{2_+}$.
    \label{lem:classspaceintersect}
\end{Remark}
\begin{Notation}
    For $\sF$ a family of subgroups of $G$, and $H$ a subgroup of $G$, denote by $\sF_H$ the family of subgroups of $H$ consisting of those groups in  $\sF$ that are contained in $H$. 
    
    If $\sF$ is a family that makes sense for all groups $G$, such as the family of all subgroups, the family of elementary abelian $p$-groups, etc., we write $\sF(G)$ for this family of subgroups of $G$. For instance, we write $\sE_{(2)}(D_8)$ for the elementary abelian subgroups of the dihedral group of order 8.
\end{Notation}

\begin{Remark}
    If $\sF$ is a family of subgroups of $G$, and $H$ is a subgroup of $G$, then $\Res^G_H E\sF \simeq E\sF_H$.
    \label{lem:restrictingclassifyingspaces}
\end{Remark}
 \begin{Lemma}
    Let $M$ be an $\sF$-nilpotent $G$-spectrum, and let $H \subset G$ be a subgroup. Then $\Res^G_H M$ is $\sF_H$-nilpotent, and $\exp_{\sF_H}\Res^G_H M \leq \exp_{\sF}M$. 
    \label{lem:expsubgp}
\end{Lemma}
\begin{proof}
    This follows from \cite[Cor.\ 4.13]{mnnnd} and \cref{lem:restrictingclassifyingspaces}. 
\end{proof}
 For a group $G$, we denote the spectrum representing Borel $G$-equivariant $\Ffield_p$-cohomology by $\underline{H\Ffield_p}_{G}$.
\begin{Corollary}
    Let $G$ be a group and let $H \subset G$ be a subgroup. Then
    \begin{equation*}
        \exp_{\sE_{(p)}(H)} \underline{H\Ffield_p}_H \leq \exp_{\sE_{(p)}(G)} \underline{H\Ffield_p}_G.
    \end{equation*}
    \label{cor:expofsubgroup}
\end{Corollary}
\begin{Lemma}
    Let $\sF_1$, $\sF_2$ be two families of subgroups, and let $M$ be a $G$-spectrum which is both $\sF_1$- and $\sF_2$-nilpotent, with exponents $m$, $n$ respectively. Then $M$ is $\sF_1 \cap \sF_2$-nilpotent, and
    \begin{equation*}
        \exp_{\sF_1 \cap \sF_2} M \leq m+n - 1.
    \end{equation*}
    \label{lem:expintersect}
\end{Lemma}
\begin{proof}
    The fact that $M$ is $\sF_1 \cap \sF_2$-nilpotent is part of \cite[Prop.\ 6.39]{mnnnd}. The assumption on the exponents implies that both maps in 
    \begin{align}
        \sk_{m-1}E\sF_{1_+} \wedge \sk_{n-1}E\sF_{2_+} \wedge M & \to \sk_{m-1}E\sF_{1_+} \wedge E\sF_{2_+} \wedge M  \label{explemmascap:eq0} \\
        & \to E\sF_{1_+} \wedge E\sF_{2_+} \wedge M
        \label{explemmascap:eq1}
    \end{align}
    have a retraction, hence the composite has a retraction. The composite of (\ref{explemmascap:eq0}) and (\ref{explemmascap:eq1}) factors as
    \begin{equation*}
        \begin{tikzcd}
            \sk_{m-1}E\sF_{1_+} \wedge \sk_{n-1}E\sF_{2_+}  \arrow{r}{\mathrm{(*)}}  \arrow{d}{\mathrm{(**)}} 
            & E\sF_{1_+} \wedge E\sF_{2_+} \wedge M \\
            \sk_{m+n-2}\left(E\sF_{1_+} \wedge E\sF_{2_+}\right) \wedge M  \arrow{ur}{\mathrm{(***)}} &
        \end{tikzcd}.
    \end{equation*}
    Composing the retraction of $\mathrm{(*)}$ with $\mathrm{(**)}$ gives a retraction of $\mathrm{(***)}$. Since $E\sF_{1_+} \wedge E\sF_{2_+} \simeq E(\sF_1 \cap \sF_2)_+$ by \cref{lem:classspaceintersect}, we obtain the desired bound. \end{proof}
\begin{Notation}
    For a $G$-space $X$, we denote by $\sI(X)$ the minimal family containing the isotropy groups of $X$. 
\end{Notation}
\begin{Example}
    For an orthogonal $G$-representation $V$, the unit sphere $S(V)$ in $V$ inherits a $G$-action. Then $\sI(S(V))$ is the smallest family containing the isotropy groups of $S(V)$.
\end{Example}
\begin{Lemma}[{{Cf. Proof of \cite[Thm.\ 2.3]{mnn}}}]
    Let $R$ be a ring $G$-spectrum with multiplicative Thom classes (e.g., $\underline{H\Ffield_2}$, see \cite[Def.\ 5.1]{mnn}). Let $V$ be a $G$-representation with corresponding oriented Euler class
    \begin{equation*}
        (\chi(V) \colon S^{-|V|} \to R) \in R^*.
    \end{equation*}
    Suppose $\chi(V)$ is nilpotent with $\chi(V)^n = 0$. Then $R$ is $\sI(S(V))$-nilpotent and 
    \begin{equation*}
      \exp_{\sI(V)}R \leq n\dim_{\Reals} V.
    \end{equation*}
   \label{cor:hfpeulerclass}
\end{Lemma}
\begin{proof}
    The fact that $\chi(V)^n = 0$ is equivalent to the oriented Euler class
    \begin{equation*}
        R \xrightarrow{e_{nV}} S^{nV} \wedge R
\end{equation*}
    being nullhomotopic (see \cite[Lem.\ 5.3]{mnn} and \cite[Rem.\ 2.4]{mnn}). Hence the left map in the fiber sequence
    \begin{equation*}
        S(nV)_+ \wedge R \to R \to S^{nV} \wedge R
    \end{equation*}
    has a section: $R$ is a retract of $S(nV)_+ \wedge R$. But $S(nV)$ is an $(n \dim_{\Reals}V - 1)$-dimensional $G$-CW complex, thus by  \cref{lem:cwexp}, $\exp_{\sF}S(nV)_+ \leq n \dim_{\Reals}V$, and hence the same bound holds for $R$ by \cref{lem:expretracts}.
   \end{proof}
\begin{Lemma}
    Let $f \colon G\to C_2 \cong O(1)$ be a 1-dimensional real representation with oriented Euler class $e \in H^1(BG;\Ffield_2)$. Suppose $e$ is nilpotent with $n$ the minimal integer $\geq 0$ such that $e^n = 0$. Then $\underline{H\Ffield_2}$ is nilpotent for the family $\All_{\ker f}$ of subgroups of $\ker f$ with $\exp_{\All_{\ker f}} \underline{H\Ffield_2} = n$.
    \label{cor:1dimrepexps}
\end{Lemma}
\begin{proof}
    The upper bound for the exponent follows immediately from \cref{cor:hfpeulerclass}.

    For the lower bound on the exponent, denote a class detecting $e$ on the $E_\infty$-page of the $\All_{\ker f}$-limit spectral sequence converging to $H^*(BG;\Ffield_2)$ by $\etilde$. Since $e$ restricts to 0 on $\ker f$, hence on all subgroups in $\All_{\ker f}$, we know that $\etilde$ lives in filtration degree $\geq 1$ on $E_\infty$. By assumption, $e^{n-1} \neq 0$. Together with the upper bound, this implies $\etilde^{n-1} \neq 0$, hence $\exp_{\All_{\ker f}} R \geq n$.
\end{proof}
We now recall the following well-known fact.
\begin{Lemma}
    For a $G$-space $X$, $H \subset G$ a subgroup, the composite of maps of Borel cohomology rings
    \begin{equation*}
        H_G^*(X; \Ffield_p) \xrightarrow{\Res^G_H} H_H^*(X;\Ffield_p) \xrightarrow{\Ind_H^G} H_G^*(X;\Ffield_p)
\end{equation*}
    is multiplication by $[G:H]$.
    \label{lem:borelresind}
\end{Lemma}

\begin{Corollary}
    Let $G$ be a group and $P \subset G$ a $p$-Sylow subgroup. Then the maps of $G$-spectra
    \begin{equation}
        \underline{H\Ffield_p}_G \to G/P_+ \wedge \underline{H\Ffield_p}_G \to \underline{H\Ffield_p}_G
        \label{eq:psylowretract}
    \end{equation}
    representing the natural transformations $\Res^G_P$ and $\Ind_P^G$ exhibit $\underline{H\Ffield_p}_G$ as a retract of $G/P_+ \wedge \underline{H\Ffield_p}_G$.
    \label{cor:psylowretract}
\end{Corollary}
\begin{proof}
    We need to show that for all subgroups $G' \subset G$, applying $\pi_*^{G'}$ to (\ref{eq:psylowretract}) yields an isomorphism. This composite is given precisely by \cref{lem:borelresind} for $X=G/G'$ and $H=P$. But multiplication by $[G:P]$ is an isomorphism because $p \nmid [G:P]$, whence $[G:P] \in \Ffield_p^{\times}$.
\end{proof}
\begin{Lemma}
    Let $G$ be a group, and write $P$ for a $p$-Sylow subgroup of $G$. Then 
    \begin{equation*}
        \exp_{\sE_{(p)}(G)} \underline{H\Ffield_p}_G = \exp_{\sE_{(p)}(P)} \underline{H\Ffield_p}_P.
    \end{equation*}
    \label{lem:psylowhfpexp}
\end{Lemma}
\begin{proof}
    First, we have
    \begin{align}
        \Res^G_{P} \underline{H\Ffield_p}_G = \underline{H\Ffield_p}_{P},  \label{explemmas:eq3} \\
        \intertext{and}
        \Res^G_{P} E\sE_{(p)}(G) = E\sE_{(p)}(P).       \label{explemmas:eq2}
    \end{align}
    Hence by \cref{cor:expofsubgroup} (cf. \cite[Cor.\ 4.13]{mnnnd}), 
    \begin{equation*}
        \exp_{\sE_{(p)}(P)} \underline{H\Ffield_p}_P \leq \exp_{\sE_{(p)}(G)} \underline{H\Ffield_p}_G.
    \end{equation*}

    For the upper bound, write $n = \exp_{\sE_{(p)}(P)} \underline{H\Ffield_p}_P$. We then have that 
    \begin{equation}
        \sk_{n-1}E \sE_{(p)}(P) \wedge \underline{H\Ffield_p}_P \to \underline{H\Ffield_p}_P
        \label{explemmas:eq1}
    \end{equation}
    admits a section.
    Note that for every $X \in \Sp_G$, we have $\Ind_H^G \Res^G_H X = G/H_+ \wedge X$. Furthermore, if we use the model for $E\sE_{(p)}(P)$ from (\ref{explemmas:eq2}), we have
    \begin{equation*}
        \sk_{n-1}E\sE_{(p)}(P) = \Res^G_P\sk_{n-1}E\sE_{(p)}(G).
\end{equation*}
    Applying this, the fact that $\Res^G_P$ is monoidal, and (\ref{explemmas:eq3}) to (\ref{explemmas:eq1}) yields a section of
    \begin{equation*}
        G/P_+ \wedge \sk_{n-1}E \sE_{(p)}(G) \wedge \underline{H\Ffield_p}_G \to G/P_+ \wedge \underline{H\Ffield_p}_G.
        \end{equation*}
Hence
    \begin{equation*}
        \exp_{\sE_{(p)}(G)} G/P_+ \wedge \underline{H\Ffield_p}_G \leq n.
    \end{equation*}
    Applying \cref{cor:psylowretract} and \cref{lem:expretracts} yields
    \begin{equation*}
        \exp_{\sE_{(p)}} \underline{H\Ffield_p}_G \leq n.
    \end{equation*}
\end{proof}
\begin{Lemma}
\label{lem:twofamexp}
    Let $\sF$ and $\sG$ be families of subgroups, and let $M$ be both $\sF$- and $\sG$-nilpotent. Then for every $K \in \sF$, $\Res_K^G M$ is $\sG_K$-nilpotent by \cref{lem:expsubgp}. Write $n = \exp_{\sF} M$, $m_K = \exp_{\sG_K} \Res^G_K M$ for all $K \in \sF$, and $m = \max_K m_K$. Then
    \begin{equation*}
        \exp_{\sG}M \leq mn.
    \end{equation*}
\end{Lemma}
\begin{proof}
    The $\sG_K$-nilpotence of $\Res^G_K M$ implies that, for all $K \in \sF$,
    \begin{equation*}
        \sk_{m_K-1}E\sG_K \wedge \Res^G_K M \to \Res^G_K M
    \end{equation*}
    admits a section. By inducing both sides from $K$ to $G$, we see that $\exp_{\sG}(G/K_+ \wedge M) \leq m_K$. By taking the coproduct over all $K \in \sF$, we obtain that $\exp_{\sG} \sk_0 E\sF \wedge M \leq m$. We now argue by induction that
    \begin{equation}
        \exp_{\sG} \sk_{d-1} E\sF \wedge M \leq md,
        \label{eq:FGind}
    \end{equation}
    for all $d \geq 1$. Hence assume (\ref{eq:FGind}) has been established for some $d \geq 1$, and consider the cofiber sequence
    \begin{equation*}
        \sk_{d-1}E\sF\wedge M \to \sk_d E\sF\wedge M \to \bigvee(S^d \wedge G/K_+)\wedge M
    \end{equation*}
    that ends in a wedge of spheres with isotropy in $\sF$ smashed with $M$. By the induction hypothesis, the $\sG$-exponent of the left term is $\leq md$, and the $\sG$-exponent of the right hand side is $\leq m$ because we already saw that $\exp_{\sG}(G/K_+ \wedge M) \leq m$ for all $K \in \sF$. Hence by \cite[Prop.\ 4.9.2]{mnnnd} the $\sG$-exponent of the middle term is $\leq m(d+1) $, which completes the induction.
    
    We have by $\sF$-nilpotency of $M$ that $M$ is a retract of $\sk_{n-1}E\sF \wedge M$, hence taking $d=n$ in (\ref{eq:FGind}) yields the result.
\end{proof}
          \subsection{Representations and exponents}
        We will now present \cref{projbund:prop1}, which will be the main tool in establishing \cref{thm:A} and \cref{thm:B}. This theorem is based on the argument in \cite[Ex.\ 5.18]{mnn} for giving upper bounds on the $\sE_{(2)}$-exponent on $\underline{H\Ffield_2}_{Q_8}$ using the projective bundle theorem. We will repeatedly use \cref{projbund:prop1} to prove upper bounds on the $\sF$ exponent of $\underline{H\Ffield_2}_G$ for various 2-groups $G$ and families $\sF \supset \sE_{(2)}$. We will also discuss a complex analogue of \cref{projbund:prop1}, which is given by \cref{cpxprojbund:prop1}.
\begin{Definition}
    Let $V$ be a real or complex vector space. Then the \emph{projectivation} $\Proj(V)$ of $V$ is the space of all respectively real or complex lines in $V$.
\end{Definition}
\begin{Remark}
    Note that if $V$ comes equipped with a linear $G$-action (i.e., is a $G$-representation), then $\Proj(V)$ inherits a natural $G$-action.
\end{Remark}
The goal of this subsection is to prove:
\begin{Theorem}
    Let $G$ be a finite group, and let $n \geq 0$ be an integer. Suppose $G$ has a real $n$-dimensional representation $V$ such that the projectivation $\Proj(V)$ has isotropy groups contained in some family $\sF$, that is, for every real line $L \subset V$ the isotropy group $G_L \leq G$ of elements of $G$ fixing $L$ satisfies $G_L \in \sF$. Then $\underline{H\Ffield_2}$ is $\sF$-nilpotent and the exponent satisfies $\exp_{\sF} \underline{H\Ffield_2} \leq n$.
    \label{projbund:prop1}
\end{Theorem}
For complex bundles we have an analogous result for cohomology with $\Integers$-coefficients:
\begin{Theorem}
    Let $G$ be a finite group, and let $n \geq 0$ be an integer. Suppose $G$ has a complex $n$-dimensional representation $V$ such that the projectivation $\Proj(V)$ has isotropy groups contained in some family $\sF$, that is, for every complex line $L \subset V$ the isotropy group $G_L \leq G$ of elements of $G$ fixing $L$ satisfies $G_L \in \sF$. Then $\underline{H\Integers}$ is $\sF$-nilpotent, and the exponent satisfies $\exp_{\sF} \underline{H\Integers} \leq 2n -1$.
    \label{cpxprojbund:prop1}
\end{Theorem}
 \subsubsection{The projective bundle theorem}
To prove \cref{projbund:prop1} and \cref{cpxprojbund:prop1}, we will use the projective bundle theorem, which can be used to develop the theory of Stiefel--Whitney classes and Chern classes. This is carried out, for instance, in \cite[Ch. 17]{husemoeller1994fibre}. We will follow the treatment and notation of \cite[\S 17.2]{husemoeller1994fibre}, but only discuss what we need. 

As in \cite[Ch.~17]{husemoeller1994fibre}, we consider real and complex vector bundles at the same time. For the case of real vector bundles, we write $c = 1$, we consider cohomology with coefficients in $K_c = \Integers/2$, and we let $F$ be the field $\Reals$ of real numbers. For the case of complex vector bundles, we write $c=2$,  we consider cohomology with coefficients in $K_c = \Integers$, and we let $F$ be the field $\Complex$ of complex numbers.

We will write $E(\eta)$ (resp.\ $B(\eta)$) for the total (resp.\ base) space, of a fiber (not necessarily vector) bundle $\eta$. Let $\xi \colon E \xrightarrow{p} B$ be an $n$-dimensional vector bundle. Let $E_{0}$ be the non-zero vectors in $E$. Let $\Proj(E)$ be the quotient of $E_0$ where we identify non-zero vectors in the same line. This yields a factorization
\begin{equation*}
    \begin{tikzcd}
E_0 \arrow{r} & \Proj(E) \arrow["q"]{r} & B,
    \end{tikzcd}
\end{equation*}
and $\Proj(E) \xrightarrow{q} B$ is a fiber bundle with fiber $F\Proj^{n-1}$, called the projectivation of $\xi$ and denoted $\Proj(\xi)$. This space admits a canonical line bundle, classified by a map $f \colon \Proj(E(\xi)) \to F\Proj^\infty$. Pulling back a polynomial generator $z \in H^*(F\Proj^\infty; K_c)$ with $|z| = c$ along $f$ gives a class $a_\xi \coloneqq f^*(z) \in H^*(\Proj(E); K_c)$.

\begin{Theorem}[{{Projective bundle theorem, see \cite[Thm.\ 17.2.5]{husemoeller1994fibre}}}]
    \label{thm:projectivebundle}
    For an $n$-dimensional vector bundle $\xi$, the classes $1, a_\xi, \ldots, a_\xi^{n-1}$ form a basis of the free $H^*(B(\xi);K_c)$-module $H^*(\Proj(E(\xi)); K_c)$. In particular, 
    \begin{equation*} 
    q^* \colon H^*(B(\xi);K_c) \to H^*(\Proj(E(\xi)); K_c)
    \end{equation*}
    is the inclusion of a summand.
\end{Theorem}
\subsubsection{From representations to exponents}
We are now ready to prove \cref{projbund:prop1} and \cref{cpxprojbund:prop1}.
\begin{proof}[{{Proof of \cref{projbund:prop1} and \cref{cpxprojbund:prop1}}}]
    Consider the Borel construction on $V$ for an arbitrary subgroup $H \leq G$, and call the resulting bundle $\xi_H$:
    \begin{equation*}
        \xi_H \colon V \to V \times_H EG\to BH.
    \end{equation*}
    Note that this is natural with respect to inclusions of subgroups. The associated projective bundle is 
    \begin{equation*}
        \Proj(\xi_H) \colon \Proj(V) \to \Proj(V \times_H EG) \to BH.
    \end{equation*}
    Observe that $\Proj(V \times_H EG) = \Proj(V) \times_H EG$. Hence by \cref{thm:projectivebundle}, 
    \begin{equation*}
    F((\Proj(V) \times_H EG)_+, HK_c) \cong F(\Proj(V)_+, \underline{HK_c})^H
\end{equation*}
is a free $F(BH,K_c) = \underline{HK_c}^H$-module. Recall that the \emph{real} dimension of $V$ is $cn$. A basis is given by the proof of \cite[Prop.\ 17.3.3]{husemoeller1994fibre}, which shows that there are classes $1, a_{\xi_H}, \ldots, a_{\xi_H}^{n-1}$ that form a basis of $\pi_*^H F(\Proj(V),\underline{HK_c})$ as a $\pi_*^H \underline{HK_c}$-module. Moreover, the element $a_{\xi_H}$ is natural with respect to restriction to subgroups, because the canonical line bundle on $\Proj(V \times_G EG)$ pulls back to the canonical line bundle on $\Proj(V \times_H EG)$.  In particular,
    \begin{equation*}
        \Res^G_H a_{\xi_G} = a_{\xi_H},
    \end{equation*}
    for all $H \leq G$. It follows that all basis elements, being powers of $a_{\xi_H}$, are natural with respect to restriction to subgroups. Hence $F(\Proj(V), \underline{HK_c})$ is a free $\underline{HK_c}$-module, because freeness of modules is an algebraic condition. Therefore, a suspension of $\underline{HK_c}$ is a retract of $F(\Proj(V),\underline{HK_c})$. But by \cref{expprop} (6), the $\sF$-exponent is invariant under suspension.     
    Therefore, by \cite[Prop.\ 4.9]{mnnnd},\begin{equation*}
        \exp_{\sF}\underline{HK_c} \leq \exp_{\sF}F(\Proj(V),\underline{HK_c}).
\end{equation*}
Furthermore, $V$ is $cn$-dimensional, hence $\Proj(V)$ admits the structure of an $(cn-c)$-dimensional $G$-CW-complex \cite[Cor.~7.2]{illman83} with isotropy contained in $\sF$ by assumption, which implies by \cref{expprop} (1) and (4) that
    \begin{align*}
        \exp_{\sF} \underline{HK_c} & \leq cn - c +1 \\
         &  = \begin{cases}
             n & \text{ if } c = 1, \\
             2n - 1& \text{ if } c = 2.
         \end{cases}
    \end{align*}
\end{proof}
         \subsection{Proper subgroups}
        Let $G$ be a finite non-abelian $2$-group of order $2^k$. The goal of this subsection is to prove
\begin{equation*}
    \exp_{\sP} \underline{H\Ffield_2}_G \leq 2\floor{\sqrt{|G| - 1}}-1,
\end{equation*}
see \cref{cor:propsubgpexp} below. This result will not be used for establishing the main results of this paper.

The argument is an adaption of the ones found in \cite[Lem.\ 4.3]{pakianathan03} and \cite{symonds91}.
\begin{Lemma}
    Let $G$ be a finite non-abelian $p$-group of order $p^k$. Then $G$ has an irreducible complex representation $V$ with $\dim_{\Complex} V \geq 2$, and moreover, all such $V$ satisfy
    \begin{equation*}
        \dim_{\Complex} V \leq \floor{\sqrt{|G| - 1}}.
    \end{equation*}
    \label{lem:irredrepspgps}
\end{Lemma}
\begin{proof}
Denote by $n_1,\ldots,n_h$ the $\Complex$-dimensions of the irreducible $\Complex$-representations of $G$. We have
    \begin{equation*}
        n_1^2 + n_2^2 + \cdots + n_h^2 = |G|
    \end{equation*}
    (see, e.g., \cite[Cor.\ 2.4.2]{Serrelinrep}). Since $G$ is non-abelian, there exist $i,j$ such that $n_i \geq 2$ \cite[Thm.\ 3.1.9]{Serrelinrep} and $n_j = 1$, corresponding to the trivial representation. Assume without loss of generality that they are respectively $n_1$ and $n_2$. Then
    \begin{align*}
        |G| & = n_1^2 + n_2^2 +  \cdots + n_h^2 \\
        & \geq n_1^2 + 1,
    \end{align*}
    hence
    \begin{equation*}
        n_1 \leq \sqrt{|G| - 1},
    \end{equation*}
    where this upper bound is in general far from optimal. Applying the floor function on both sides preserves the inequality and does not change the integer on the left hand side. This yields the result.
\end{proof}
\begin{Corollary}
    Let $G$ be a finite non-abelian $2$-group of order $2^k$, and let $\sP$ be the family of proper subgroups of $G$. Then
    \begin{equation*}
        \exp_{\sP} \underline{H\Ffield_2}_G \leq 2 \floor{\sqrt{|G|-1}} -1.
    \end{equation*}
    \label{cor:propsubgpexp}
\end{Corollary}
\begin{proof}
    Let $V$ be an irreducible complex representation of $G$ satisfying
    \begin{equation*}
        2 \leq \dim_{\Complex} V \leq \floor{\sqrt{|G| - 1}},
    \end{equation*}
    as furnished by \cref{lem:irredrepspgps}. Then the complex projectivation $\Proj(V)$ has isotropy in $\sP$ the family of proper subgroups of $G$, for if $L \in \Proj(V)$ were fixed by all of $G$, $V$ would not be irreducible, since $\dim_{\Complex} V \geq 2$. An application of \cref{cpxprojbund:prop1} yields the result.
\end{proof}

     \section{Exponents of small 2-groups}
        We now prove \cref{thm:B} from the introduction, by proving the claimed upper bounds of the $\sE_{(2)}$-exponents. We first treat the class of abelian groups, dihedral groups, and the generalized quaternion groups. The remaining groups of order 16 are treated using a case-by-case analysis.
         \label{sec:expsm2gps}
        \subsection{Abelian 2-groups}
        We determine an upper bound on the $\sE_{(2)}$-exponent of $\underline{H\Ffield_2}_A$ for $A$ a finite abelian group.

\begin{Proposition}
    \label{prop:expelab}
    Let $A$ be a group isomorphic to
    \begin{equation*}
      \prod_{j \in J} C_{2^{n_j}}  \times \prod_{k \in K} C_2
    \end{equation*}
    with $n_j \geq 2$.
    Then
    \begin{equation*}
        \exp_{\sE_{(2)}}\underline{H\Ffield_2}_A \leq  \#J + 1.
\end{equation*}
\end{Proposition}
\begin{proof}
    For $j \in  J$, consider the projection of $A$ onto the $C_2$ in the $j$-th factor:
\begin{equation*}
    p_j \colon A \to C_2 \cong O(1).
\end{equation*}
The corresponding Euler class is $a_j \in H^1(BA; \Ffield_2)$. Since $a_j^2 = 0$, we obtain 
\begin{equation*}
\exp_{\All_{\ker p_j}} \underline{H\Ffield_2} = 2.
\end{equation*}
Now $\bigcap_{j  \in J}\ker p_j = \sE_{(2)}$, and hence by \cref{lem:expintersect}, we obtain $\exp_{\sE_{(2)}} \underline{H\Ffield_2} \leq \# J+1$. 
\end{proof}
         \subsection{Dihedral groups}
        Let $D_{2^n}$ be the dihedral group of order $2^n$ $(n \geq 3)$ with presentation
\begin{equation*}
    D_{2^n} = \langle \sigma, \, \rho \, | \, \sigma^2 = \rho^{2^{n-1}} = e, \, \sigma \rho \sigma^{-1} = \rho^{-1} \rangle.
\end{equation*}
Denote by $T$ the matrix representing reflection in the $x$-axis. In the remainder of this section, for an angle $\theta$, denote by $R_\theta$ the matrix representing counterclockwise rotation in $\Reals^2$ by $\theta$ about the origin.
Then $\sigma \mapsto T$,  $\rho \mapsto R_{2\pi/2^{n-1}}$, yields a linear 2-dimensional real orthogonal representation of $D_{2^n}$, which we call $V$.  We apply \cref{projbund:prop1} to compute an upper bound on the exponent with respect to the family $\sE_{(2)}$ of elementary abelian 2-subgroups. In the following lemma we establish that the isotropy of $\Proj(V)$ is contained in $\sE_{(2)}$. 
\begin{Lemma}
    The isotropy groups of the projectivation $\Proj(V)$ are contained in the family $\sE_{(2)}$ of elementary abelian 2-groups in $D_{2^n}$.
    \label{lem:d8isotropy}
\end{Lemma}
\begin{proof}
    The proof is by elementary linear algebra.
    Since $V$ is an orthogonal representation, a linear subspace $L$ spanned by a vector $v$ is fixed by an element $g$ of $D_{2^n}$ if and only if $g$ has eigenvalue $1$ or $-1$.
    A finite group is an elementary abelian 2-group if and only if all of its elements have order dividing 2. Therefore, we can reduce to proving that all elements of order 4 do not fix any lines, which we can reduce further to proving that all elements of order 4, up to inverses and conjugation, do not fix any lines. This means it suffices to check that the element $\rho^{2^{n-3}}$ does not fix any lines. But this element has characteristic polynomial $\lambda^2 + 1$, which has no real zeroes, and thus the result follows.
\end{proof}
\begin{Corollary}
    The $\sE_{(2)}$-exponent of $\underline{H\Ffield_2}_{D_{2^n}}$ satisfies
    \begin{equation*}
        \exp_{\sE_{(2)}} \underline{H\Ffield_2}_{D_{2^n}} \leq 2.
    \end{equation*}
    \label{cor:expd8upperbound}
\end{Corollary}
\begin{proof}
    Immediate from \cref{projbund:prop1} and \cref{lem:d8isotropy} and the fact that $V$ is 2-dimensional.
\end{proof}

         \subsection{Quaternion groups}
        The generalized quaternion group $Q_{2^n}$ of order $2^n$ is the finite subgroup of quaternionic space $\HH$ generated multiplicatively by the elements of unit length $\{ e^{2\pi i / 2^{n-1}},\, j\}$ (see, e.g., \cite[XII.\S7]{cartaneilenberg}). Denoting these generators by $r$ and $s$, respectively, one obtains the presentation
\begin{equation*}
    Q_{2^n} = \langle r,\, s\, | \, r^{2^{n-2}} = s^2,\, rsr = s \rangle,
\end{equation*}
for all $n \geq 3$. The subgroup generated by $\langle s^2 \rangle$ is central and of order 2, which gives rise to a central extension \cite[IV.2]{adem}
\begin{equation*}
    C_2 \to Q_{2^n} \to D_{2^{n-1}}.
\end{equation*}
In this subsection we will prove the following upper bound on the $\sE_{(2)}$-exponent.
\begin{Proposition}
    The $\sE_{(2)}$-exponent satisfies
    \begin{equation*}
        \exp_{\sE_{(2)}} \underline{H\Ffield_2}_{Q_{2^n}} \leq 4.
    \end{equation*}
    \label{q2n:expupper}
\end{Proposition}
\begin{proof}
    The proof is a straightforward adaption of the argument in \cite[Ex.\ 5.18]{mnn}.

    Let $\HH \cong \Reals^4$ be the 4-dimensional real representation coming from the embedding $Q_{2^n} \hookrightarrow \HH$. This is a free action, and restricts to a free action on $S^3$. The subgroup $\langle \pm 1 \rangle$ is central, and therefore $Q_{2^n}/\langle \pm 1 \rangle$ acts on $S^3/\langle \pm 1 \rangle =  \Proj(\Reals^4)$ with isotropy in $\langle \pm 1 \rangle$. The result now follows from \cref{projbund:prop1}.
\end{proof}
         \subsection{The modular group of order 16}
        Let
\begin{equation*}
    M_{16}= \langle r,\, f \, | \, r^8 = f^2 = e, \, frf^{-1}=r^5 \rangle
\end{equation*}
be the modular group of order 16. The group has this name because its lattice of subgroups is modular (see, e.g., \cite[I.\S7]{birkhoff67}).

In this subsection we prove the following upper bound on the $\sE_{(2)}$-exponent.
\begin{Proposition}
    The $\sE_{(2)}$-exponent satisfies
    \begin{equation*}
        \exp_{\sE_{(2)}} \underline{H\Ffield_2}_{M_{16}} \leq 4.
    \end{equation*}
    \label{m16:expupper}
\end{Proposition}
\begin{proof}
    For an angle $\theta$, let $R_{\theta}$ be the rotation matrix. Let $T$ be the matrix 
which interchanges the summands of $\Reals^2 \oplus \Reals^2$. Write $\alpha = 2\pi/8$. Let $V$ be the 4-dimensional real orthogonal $M_{16}$-representation given by $f \mapsto T$, $r \mapsto R_\alpha \oplus R_{5\alpha}$.
This is the representation from \cite[Lem.\ 13.3]{totaro14}.

    We will now show that the isotropy groups of the projectivation $\Proj(V)$ are contained in $\sE_{(2)}$. As in the proof of \cref{lem:d8isotropy}, we only need to consider all the elements of order 4, up to taking inverses and conjugation, of which there are two: $r^2$ and $fr^2$. Both have characteristic polynomial $(\lambda^2 + 1)^2$, which has no real roots. Therefore, they do not fix any lines.
Applying \cref{projbund:prop1} gives the desired result.
\end{proof}

         \subsection{The semi-dihedral group of order 16}
        Let
\begin{equation*}
    SD_{16} = \langle s,r \, | \, s^2 = r^8 = e, \,  srs^{-1} = r^3 \rangle
\end{equation*}
be the semidihedral group of order 16. 
\begin{Proposition}
    The $\sE_{(2)}$-exponent satisfies
    \begin{equation*}
        \exp_{\sE_{(2)}} \underline{H\Ffield_2}_{SD_{16}} \leq 4.
    \end{equation*}
    \label{prop:sd16upper}
\end{Proposition}
\begin{proof}
    We will use the notation $R_\theta$ and $T$ from the proof of \cref{m16:expupper}. 
Write $\alpha = 2 \pi /8$ and let $V$ be the real orthogonal $SD_{16}$-representation given by $s \mapsto T$, $r \mapsto R_\alpha \oplus R_\alpha$
from \cite[Lem.\ 13.4]{totaro14}. 
    
    Again, it reduces to considering two elements: $r^2$ and $sr$. The element $r^2$ has characteristic polynomial $(\lambda^2 + 1)^2$, and the element $sr$ has characteristic polynomial $\lambda^4 + 1$. Neither of these has real roots, so these elements do not fix any lines.
Applying \cref{projbund:prop1} gives the desired result.
\end{proof}

         \subsection{The central product of $D_8$ and $C_4$}
        Let $D_8 = \langle \sigma, \rho \rangle$ be the dihedral group of order 8, and let $C_4 = \langle \gamma \rangle$ be the cyclic group of order 4. Both these groups have central cyclic subgroups of order 2: for $D_8$ this is $\langle \rho^2 \rangle$ and for $C_4$ this is $\langle \gamma ^2 \rangle$. The \emph{central product} of $D_8$ and $C_4$ is defined to be the direct product with these central subgroups identified:
\begin{equation*}
    D_8 \ast C_4 := D_8 \times C_4 / \langle \rho^2 \gamma^{-2} \rangle.
\end{equation*}
\begin{Proposition}
    The $\sE_{(2)}$-exponent satisfies $\exp_{\sE_{(2)}} \underline{H\Ffield_2}_{D_8 \ast C_4} \leq 4$.
    \label{d8cpc4:elexpupperbound2}
\end{Proposition}
\begin{proof}
    Let $R_{\pi/2}$ be the rotation matrix with $\theta = \pi/2$, and let $T$ be the matrix interchanging the summands of $\Reals^2 \oplus \Reals^2$. Consider the representation given by $\sigma \mapsto \Id \oplus(-\Id)$, $\rho \mapsto T((-\Id) \oplus \Id)$, $\gamma \mapsto R_{\pi/2} \oplus R_{\pi/2}$,
which is the underlying 4-dimensional real representation from \cite[Lem.\ 13.5]{totaro14}. 

Again, it suffices to study the elements $\rho$, $\gamma$, and $\rho^2 \gamma$. The element $\rho$ has characteristic polynomial $\lambda^4 + 1$, the element $\gamma$ has characteristic polynomial $(\lambda^2 + 1)^2$. Neither of these polynomials has real roots, so neither $\rho$ nor $\gamma$ fixes any lines. Since $\rho^2$ acts by $-\Id$ it follows that also $\rho^2 \gamma$ does not fix any lines.
The result follows from \cref{projbund:prop1}.
\end{proof}
         \subsection{The group $(C_4 \times C_2) \rtimes C_2$}
        Let $C_4 \cong \langle a \rangle$, $C_2 \cong \langle b \rangle$, and let $C_2 \cong \langle c \rangle$ act on $C_4 \times C_2 \cong \langle a, b \rangle$ via $a \mapsto ab$, $b \mapsto b$.
We consider the group $G$ which is defined to be the semidirect product of $(C_4 \times C_2) \rtimes C_2$ under this action.
This group hence admits the following presentation:
\begin{equation*}
	G = \langle a,b,c \, | \, a^4 = b^2 = c^2 = e, \, ab = ba, \, bc = cb, \, cac^{-1} = ab \rangle.
\end{equation*}

We will establish an upper bound on the $\sE_{(2)}$-exponent of $G$. \begin{Proposition}
    The $\sE_{(2)}$-exponent of $G = (C_4 \times C_2) \rtimes C_2$ satisfies
    \label{prop:c4xc2sdc2exp}
\begin{equation*}
\exp_{\sE_{(2)}} \underline{H\Ffield_2} \leq 2.
\end{equation*}
\end{Proposition}
\begin{proof}
    Let $A$ be the subgroup $\langle a^2, bc, c\rangle$, and let $B$ be the subgroup of $A$ given by $\langle bc, c\rangle$. Both $A$ and $B$ are normal subgroups of $G$. Furthermore, $A$ is the unique maximal elementary abelian 2-subgroup of $G$. Consider the composition of quotient maps
    \begin{equation*}
        G \to G/B \cong C_4 \to G/A \cong C_2.
    \end{equation*}
    Pulling back the sign representation of $C_2$ to $G$ gives a 1-dimensional real $G$-representation with Euler class $e$. Because of the factorization over $C_4$, we have $e^2 = 0$, by naturality of the Euler class, and the fact that every class in $H^1(BC_4; \Ffield_2)$ squares to 0.
Hence by \cref{cor:1dimrepexps}, $\exp_{\sE_{(2)}}\underline{H\Ffield_2} \leq 2$.
\end{proof}

         \subsection{The group $C_4 \rtimes C_4$}
        Let $r$ and $s$ generate two copies of $C_4$: $\langle r \, | \, r^4 = e \rangle,\, \langle s \, | \, s^4 = e \rangle \cong C_4$. Let $\langle s \rangle$ act on $\langle r \rangle$ by $s \cdot r = r^{-1}$. Then $C_4 \rtimes C_4$ is defined to be the corresponding semi-direct product $\langle r \rangle \rtimes \langle s \rangle$. A presentation is $C_4 \rtimes C_4 = \langle r,s \, | \, r^4 = s^4 = e, srs^{-1} = r^3 \rangle$. We will determine an upper bound on the $\sE_{(2)}$-exponent of $\underline{H\Ffield_2}_{C_4 \rtimes C_4}$. In order to do so, we will use the fact that $H^*(BC_4 \rtimes C_4;\Ffield_2)$ is isomorphic to \cite[App.\ C, \#10(16)]{carlson2003cohomology}
\begin{equation}
    \Ffield_2[z, y, x, w]/(z^2 + y^2, zy)
    \label{eq:hc4rtimesc4}
\end{equation}
with degrees $|z| = |y| = 1$, $|x| = |w| = 2$.
We will determine an upper bound on $\exp_{\sE_{(2)}} \underline{H\Ffield_2}$ using the following lemmas.
\begin{Lemma}
    The Euler class $e$ of the 1-dimensional real representation given by pulling back the sign representation along the quotient map
    \begin{equation}
        C_4 \rtimes C_4 \to C_4 \rtimes C_4/\langle r, s^2 \rangle \cong C_2
        \label{eq:c4rtc4rs2quot}
    \end{equation}
    satisfies $e^2 = 0$.
    \label{c4sdpc4:rep1}
\end{Lemma}
\begin{proof}
    The proof is the same as the proof of \cref{prop:c4xc2sdc2exp}: the quotient map (\ref{eq:c4rtc4rs2quot}) factors through the map $C_4 \rtimes C_4 \to C_4 \rtimes C_4/\langle r \rangle \cong C_4$, hence $e^2 = 0$ by naturality.
\end{proof}
\begin{Lemma}
    The Euler class $e$ of the 1-dimensional real representation given by pulling back the sign representation along the quotient map
    \begin{equation*}
        C_4 \rtimes C_4 \to C_4 \rtimes C_4/\langle r^2, s \rangle \cong C_2
    \end{equation*}
    satisfies $e^3 = 0$.
    \label{c4sdpc4:rep2}
\end{Lemma}
\begin{proof}
    This follows from the fact that that $H^*(BC_4 \rtimes C_4; \Ffield_2)$ is isomorphic to (\ref{eq:hc4rtimesc4}), and that in that graded ring every element in degree 1 cubes to 0.
\end{proof}
Combining these lemmas, we obtain
\begin{Proposition}
    \label{c4sdpc4:goodupper}
    The $\sE_{(2)}$-exponent satisfies
    \begin{equation*}
        \underline{H\Ffield_2}_{C_4 \rtimes C_4} \leq 4.
    \end{equation*}
\end{Proposition}
\begin{proof}
    We remark that $\langle r^2, s^2 \rangle$, which is the unique maximal elementary abelian 2-subgroup of $C_4 \rtimes C_4$, is the intersection of $\langle r^2, s \rangle$ and $\langle r, s^2 \rangle$, the groups we divided out by in \cref{c4sdpc4:rep1} and \cref{c4sdpc4:rep2}. Hence by \cref{cor:1dimrepexps} and \cref{lem:expintersect}, $\exp_{\sE_{(2)}} \underline{H\Ffield_2} \leq 2 + 3 - 1 = 4$.
\end{proof}
         \subsection{The group $D_8 \times C_2$}
        \begin{Proposition}
    \label{prop:d8xc2e2exp}
    For $D_8 \times C_2$, the $\sE_{(2)}$-exponent of $\underline{H\Ffield_2}$ satisfies
    \begin{equation*}
        \exp_{\sE_{(2)}} \underline{H\Ffield_2} \leq 2.
    \end{equation*}
\end{Proposition}
\begin{proof}
    Pulling back transitive $D_8$-orbits with isotropy in elementary abelian groups along the projection map $D_8 \times C_2 \to D_8$ gives transitive $D_8 \times C_2$-orbits with isotropy in elementary abelian groups. Therefore, pulling back the representation of $D_8$ considered in the proof of \cref{cor:expd8upperbound} along the projection map $D_8 \times C_2 \to D_8$ gives a 2-dimensional real representation of $D_8 \times C_2$ with projectivation having isotropy contained in $\sE_{(2)}$. Hence by \cref{projbund:prop1}, we have $\exp_{\sE_{(2)}} \underline{H\Ffield_2}_{D_8 \times C_2} \leq 2$.
\end{proof}
         \subsection{The group $Q_8 \times C_2$}
        \begin{Proposition}
    \label{prop:q8xc2e2exp}
    For $Q_8 \times C_2$, the $\sE_{(2)}$-exponent of $\underline{H\Ffield_2}$ is
    \begin{equation*}
        \exp_{\sE_{(2)}} \underline{H\Ffield_2}_{Q_8 \times C_2} \leq 4.
    \end{equation*}
\end{Proposition}
\begin{proof}
    Elementary abelian subgroups of $Q_8$ pull back to elementary abelian subgroups of $Q_8 \times C_2$ along the projection map $Q_8 \times C_2 \to Q_8$. Therefore, pulling back the 4-dimensional real representation of $Q_8$ whose projectivation has isotropy in $\sE_{(2)}$ considered in \cite[Ex.\ 5.18]{mnn} along the projection map gives a 4-dimensional real representation of $Q_8 \times C_2$ with projectivation with isotropy in $\sE_{(2)}(Q_8 \times C_2)$. Hence, by \cref{projbund:prop1}, the exponent satisfies $\exp_{\sE_{(2)}} \underline{H\Ffield_2}_{Q_8 \times C_2} \leq 4$. 
\end{proof}
     \bibliographystyle{plain}

\end{document}